\documentclass[10pt]{amsart}
\usepackage{amsmath, latexsym, 
}
\usepackage{amssymb,amsthm}
\usepackage{amscd,graphics}

\theoremstyle{remark}

\theoremstyle{definition}

\def\id{\mathrm{id}}

\def\EE{\mathcal{E}}
\def\FF{\mathcal{F}}

\def\JJ{\mathcal{J}}

\def \PP{\mathcal {P}}

\begin{document}
\title[Corrigendum]{Corrigendum to ``Linear response formula
for piecewise expanding unimodal maps," Nonlinearity, 21 (2008) 677--711}
\author{V. Baladi and D. Smania}

\date{May 24, 2012}
\maketitle


\section{Theorem 7.1}

The last claim of Theorem 7.1 should read:
``If the postcritical orbit of $f_0$ is dense in $[c_2, c_1]$ then there exist
$\varphi \in C^\infty(I)$  and a sequence
$ t_n \to 0$ so that $c$ is not periodic under any $f_{t_n}$, and
so that we have
$$
\lim_{n \to \infty}| t_n^{-1}
(\int \varphi \rho_{t_n} - \int \varphi \rho_0\,dx )|
\to \infty \, ."
$$

\begin{proof}[Proof if the orbit of $c$ is  dense]
We have  $\EE_0=\FF_0=\id$ and, using (94),
we can consider $\PP_t$ as in the case when the orbit is infinite
but not dense.

Let $x_0$ be the fixed point of $f$ which lies in the
interior of $I$, and assume that $\varphi(x_0)=1$
and $\int \varphi \, d\mu=0$.
Since the postcritical orbit is dense,  for any  $\delta > 0$ there
exists $j_0(\delta) \ge 1$ so that $d(c_{j_0}, x_0)< \delta$. Clearly,
$\lim_{\delta \to 0} j_0=\infty$.
 Put
$\Lambda_f = \sup |f'|$. If  $\delta \Lambda_f^m\le 1/2$ 
for some large $m$  then for all $j_0 \le n \le j_0+m$ we have
$$
\sum_{k=0}^{n-j_0}  |c_{j_0+k}-f^{k}(x_0)|\le
\delta \sum_{k=0}^{n-j_0} \Lambda_f^{k}
\le \frac{\delta \Lambda_f^{n-j_0}}{1-1/\Lambda_f}\, ,
$$
and thus
\begin{align*}
\nonumber |\sum_{k=1}^{n} \varphi(c_{k})|& \ge |\sum_{k=0}^{n-j_0} \varphi(f^k(x_0))|- 
 \sum_{k=j_0}^{n}|  \varphi(c_{k})-\varphi(f^{k-j_0}(x_0))| -|\sum_{k=1}^{j_0-1} \varphi(c_{k})|\\
&\ge (n-j_0+1)- 
\frac{1 }{2 (1-1/\Lambda_f)}\sup |\varphi'| -|\sum_{k=1}^{j_0-1} \varphi(c_{k})|\, .
\qquad (\star)
\end{align*}

Let $t_n\to 0$ be a sequence of non periodic parameters
and let $M(t_n)$ be defined by (92). 
Now, 
 \begin{align*}
 \sum_{ k =1}^{ M(t_n)} \varphi(c_k)\sum_{j=1}^k 
\frac{X(c_j)}{(f^{j-1})'(c_1)} 
&= \JJ(f,X) \sum_{ k =1}^{M(t_n)} \varphi(c_k) 
-
 \sum_{ k =1}^{M(t_n)} \varphi(c_k)\sum_{j=k+1}^\infty 
\frac{X(c_j)}{(f^{j-1})'(c_1)} 
\end{align*}
and
\begin{align*}
 |\sum_{ k =1}^{M(t_n)} \varphi(c_k)\sum_{j=k+1}^\infty 
\frac{X(c_j)}{(f^{j-1})'(c_1)} |
&\le \sup |X| \sup |\varphi|  
\sum_{ k =1}^{M(t_n)} \frac{(\inf |f'| )^{-k}}{1-1/\inf|f'|}\\
&\le \sup |X| \sup |\varphi|  
\frac{(\inf |f'|) ^{-1}}{(1-1/\inf|f'|)^2}
\end{align*}

Thus, for arbitrarily large $n$, 
recalling also $|C_n|\le \widehat C$ from the previous cases,
\begin{align*}
\nonumber &\left |\sum_{ k =1}^{ M(t_n)} \frac{s_{1,t_n}}{(f^{k-1}_{t_n})'(c_1)}
\int \varphi \frac{H_{c_{k,t_n}}- H_{c_k}}{t_n} dx\right |
=\left |C_{n} + s_1\sum_{ k =1}^{ M(t_n)} \varphi(c_k)\sum_{j=1}^k 
\frac{X(c_j)}{(f^{j-1})'(c_1)}\right |\\
&\ge 
|\sum_{ k=1}^{ M(t_n)} \varphi(c_k)||\JJ(f,X)|- \widehat C
- \sup|\varphi| \sup|X| 
\frac{\inf |f'| ^{-1}}{(1-1/\inf|f'|)^2}\,  .
\end{align*}

Assume now for a contradiction
that for any sequence
$t_n\to 0$ as above we have
$|\int \varphi d\mu_{t_n}|\le A |t_n|$  for some $A<\infty$ and
all large enough $n$. 
Then,  for all large enough $n$,  we would have
\begin{equation*}
|\sum_{ k=1}^{  M(t_n)} \varphi(c_k)|
\le \frac{A+\widetilde C(f_t,\varphi)}{ |\JJ(f,X)|}=:D \, .
\qquad \qquad\qquad\qquad(\star\star)
\end{equation*}
To end the proof we shall find
sequences $t_n$ so that the above estimate gives a contradiction.
(Note that $\varphi$ cannot be a coboundary since $\varphi(x_0)=1$.)

For $m\ge 1$, let  $\delta(m)>0$ be so that $\delta \Lambda^m_f<1/2$.
Next, take $j_0=j_0(\delta(m))\ge 1$ so that $d(c_{j_0}, x_0)<\delta$.
Then, letting $J(m) \ge 1$ be minimum for the property
$j_0(\delta(m+J(m)))-1  > j_0(\delta(m))+m$, we have
$$
j_0(\delta(m)) -1< j_0(\delta(m))+m < j_0(\delta(m+J(m)))-1 < \ldots\, ,
$$
and this defines a sequence, denoted $L(n)$, 
so that $L(n) \to \infty$ as $n \to \infty$.
We claim that we can choose
the sequence $t_n\to 0$ of nonperiodic parameters so that for  all large enough $n$ we have  $M(t_n)=L(n)$.  Indeed, 
by the definition of $M(t)$, and since $\inf |f'|>1$,  there is a sequence
$ 0< \tau_{L} < \tau_{L-1}$, $L \ge 1$, with $\tau_L \to 0$ as $L \to \infty$,
so that for any $t \in [\tau_L, \tau_{L-1})$ we have $M(t)=L$.
Thus, since the set of non periodic parameters is dense
(see \cite[Cor. 4.1, item A]{BS}), 
there is 
a sequence of non periodic parameters
$t_n \to 0$ so that  $M(t_n)=L(n)$.

Then, recalling ($\star$),
\begin{align*}
|\sum_{k=1}^{j_0(\delta(m))+m-1} \varphi(c_{k})|&\ge
m- \frac{1}{2(1-1/\Lambda_f)}\sup |\varphi'|- |\sum_{k=1}^{j_0(\delta(m))-1} \varphi(c_k)|\\
\nonumber &\ge m-\frac{1}{2(1-1/\Lambda_f)} \sup |\varphi'|- D \, .
\end{align*}
The rightmost lower bound  clearly diverges as $m\to \infty$, giving the
desired contradiction with ($\star\star$).
\end{proof}

\section{Typographical errors}

We use this opportunity to correct two minor typographical errors:


In the Proof of Theorem 7.1, if the orbit of $c$ is infinite but not dense, ``We now consider the first term in (95)" on p. 706, line 6, should read ``We now consider the second term in (95)."

In the beginning of $\S 3.3$, functions in $\widetilde {BV}$ are supported in 
$(-\infty,b]$, not $[a,b]$. 


 \bibliographystyle{alpha}

\end{document}